\documentclass[11pt]{amsart}
\usepackage{amsmath,amsthm,epsfig,amssymb}
\usepackage{amsthm,amsfonts,amssymb,amscd}
\usepackage[utf8x]{inputenc}
\usepackage{textcomp}
\usepackage{tikz}
\usetikzlibrary{positioning, matrix, arrows}
\usepackage[all]{xy}
\usepackage[colorlinks=true,citecolor=blue,urlcolor=blue,linkcolor=blue]{hyperref}
\usepackage{url}

\newtheorem{Theorem}[subsection]{Theorem}

\theoremstyle{definition}
\newtheorem{Remark}[subsection]{Remark}

\newtheorem{proposition-definition}[subsection]{Proposition-Definition}

\newcommand{\E}{\mathcal{F}_2(E)}
\newcommand{\Oh}{\mathcal{O}}
\newcommand{\Z}{\mathbb{Z}}

\title{Stability of Secant Bundles on Second Symmetric Power of Curves}

\author[S.Basu]{Suratno Basu}
\address{Institute of Mathematical Sciences, HBNI, CIT Campus, Tharamani, Chennai 600113,
India}
\email{suratnob@imsc.res.in}

\author[K. Dan]{Krishanu Dan}
\address{Chennai Mathematical Institute, H1 Sipcot IT Park \\
Siruseri, Kelambakkam - 603103, INDIA.}
\email{krishanud@cmi.ac.in}


\begin{document}

\begin{abstract}
Given a rank $r$ stable bundle over a smooth irreducible projective curve $C,$ 
there is an associated rank $2r$ bundle over $S^2(C),$ the second symmetric 
power of $C.$ In this article we study the slope (semi-)stability of this bundle.
\end{abstract}

\maketitle

\section{Introduction}

Let $C$ be a smooth, irreducible, projective curve of genus $g \geq 2$ over $\mathbb{C},$ 
the field of complex numbers, and let $S^n(C)$ denotes the $n$-fold symmetric 
power of $C.$ Given a vector bundle $E$ of rank $r$ on $C,$ one can associate a 
rank $nr$ vector bundle $\mathcal{F}_n(E)$ on $S^n(C).$ This bundle was first studied 
by R. Schwarzenberger (\cite{Schwarzenberger}) for the case of line bundles, called it 
{\it secant bundle}. 
He used it to study the ring of rational equivalence class of $S^n(C).$ 
It is natural to ask what properties of $E$ will be inherited by $\mathcal{F}_n(E)$. 
In particular we can ask the following questions:

\begin{enumerate}
 \item If $E$ is (semi-)stable on $C$, then does it imply 
 $\mathcal{F}_n(E)$ is (semi-)stable on $S^n(C)$ with respect to some suitably chosen 
 ample divisor?
 \item If $E$ and $F$ are two bundles on $C$ such that $\mathcal{F}_n(E)\simeq 
 \mathcal{F}_n(F)$ on $S^n(C)$, then does it imply $E\simeq F$?
\end{enumerate}

The question $(2)$ has completely been answered when $E$ and $F$ are both stable in 
\cite{BN}, \cite{BP}. 
The question $(1)$, in the case when $E$ is a line bundle, is treated in the papers 
\cite{LMN}, \cite{BN1}, \cite{M}. Also, recently in \cite{DP} the authors consider a 
rank $2$ stable bundle $E$ under some generality conditions and 
show that $\mathcal{F}_2(E)$ over $S^2(C)$ is semi-stable. However, for a general 
rank $r$ (semi-)stable bundles the question $(1)$ is not yet studied to the  
best of our knowledge. In this short note we will make an attempt to understand the 
question $(1)$ when $E$ is a (semi-)stable rank $r$ bundle and $n=2$.

More precisely, we prove the following
\begin{Theorem}\label{MT}
  Let $E$ be a rank $r$ vector bundle on $C$ of degree $d.$ 
  
  $(i)$ If $E$ is slope semi-stable and $d \geq r$, the bundle $\E$ on $S^2(C)$ is 
  slope semi-stable with respect to the ample class $x+C.$ 
 
 $(ii)$ If $E$ is slope stable and $d>r$, then the bundle $\E$ on $S^2(C)$ is slope 
 stable with respect to $x+C.$
\end{Theorem}

\section{Preliminaries}

Let $C$ be a smooth, irreducible, projective curve of genus $g \geq 2$ over the 
field of complex 
numbers $\mathbb{C}.$ On the Cartesian product $C \times C,$ we have a natural 
$\mathbb{Z}/2\mathbb{Z}$ action by means of the involution $\sigma : C \times C \to 
C \times C, (x,y) \mapsto (y,x).$ Let us denote the quotient space by $S^2(C).$ It is 
a smooth, irreducible, projective surface over $\mathbb{C},$ called the second 
symmetric power of curve. Let $\pi : C \times C \to S^2(C)$ be the quotient map. 

Let $q_1$ and $q_2$ be the projections from $S^2(C) \times C$ to $S^2(C)$ and $C$ 
respectively. Let 
$$
\Delta_2 := \{ (D, z) \in S^2(C) \times C \mid z \in \text{Supp}(D) \} \subset S^2(C) 
\times C.
$$
Then $\Delta_2$ is a smooth divisor of $S^2(C) \times C,$ called the universal effective 
divisor of degree $2$ on $C.$ Let $q : \Delta_2 \to S^2(C)$ be the restriction of $q_1.$ 
Then $\Delta_2$ is a two-sheeted ramified cover of $S^2(C)$.

Let $E$ be a rank $r$ vector bundle on $C.$ Define
$$
\mathcal{F}_2(E) := q_{*}(q_2^{*}(E)|_{\Delta_2}).
$$
Then $\mathcal{F}_2(E)$ is a rank $2r$ vector bundle on $S^2(C),$ called the secant 
bundle.

For $x \in C$, let us denote by $x+C$ the reduced divisor in $S^2(C)$ whose support 
is the set $\{x+c : c \in C\}$, and let $x$ denotes the cohomology class of $x+C$ 
in $H^2(S^2(C), \mathbb{Q})$.
Also let $\theta\in H^2(S^2(C), \mathbb{Q})$ be the pullback of the cohomology class 
corresponding to the Theta divisor in
$\text{Pic}^2(C)$ under the natural morphism $S^2(C) \to \text{Pic}^2(C)$. If 
$E$ is a rank $r$ and degree $d$ vector bundle on $C$, then the first Chern character  
of $\mathcal{F}_2(E)$ has the following expression:
\begin{equation}\label{eqn0''}
 ch(\mathcal{F}_2(E))= d(1 - e^{-x}) - r(g-1) + r(g+1+ \theta)e^{-x}
\end{equation}
(\cite [Chapter VIII, Lemma 2.5]{ACGH}, \cite[p. 774]{BL}). From this we get the first 
Chern class of $\mathcal{F}_2(E)$ as
\begin{equation}\label{eqn0}
 c_1(\mathcal{F}_2(E)) = (d - r(g+1))x + r\theta.
\end{equation}

\section{Stability of Secant Bundles}

Let $C$ be a smooth, irreducible, projective curve over $\mathbb{C}$ of genus 
$g \geq 2$ and let $E$ be a rank $r$ vector bundle of degree $d$ on $C.$ 
Throughout this section, 
(semi-)stability will always mean slope (semi-)stability.

By \cite[Chapter 3, Lemma 3.2.2]{HL}, proving the semi-stability of $\E$ on $S^2(C)$ with respect 
to $x+C$ is same as proving the semi-stability of $\pi^*(\E)$ on $C \times C$ 
with respect to the ample divisor $H:= \pi^*(x+C) = [x\times C + C \times x].$ 
Using equation (\ref{eqn0}), we can see that (\cite[p. 39]{DP})
$$
c_1(\pi^*(\mathcal{F}_2(E)))= d[x \times C + C \times x] - r\Delta
$$
where $\Delta$ is the diagonal of $C \times C$, and
\begin{equation}\label{eqn0'}
 \mu_H(\pi^*(\mathcal{F}_2(E)))=\frac{d-r}{r}.
\end{equation}

Let $p_i: C \times C \to C$ be the $i$-th co-ordinate projections, $i=1,2$. The 
vector bundles $\pi^*(\E)$ and $p_1^*(E) \oplus p_2^*(E)$ are isomorphic outside 
the diagonal $\Delta$. On $C \times C$, these two vector bundles are related by the 
following exact sequence:
\begin{equation}\label{eqn1}
 0 \rightarrow \pi^*(\E) \to p_1^*(E) \oplus p_2^*(E) 
 \xrightarrow{q} E \simeq p_1^*(E)|_{\Delta} \simeq p_2^*(E)|_{\Delta} \rightarrow 0
\end{equation}
where the homomorphism $q$ is defined as 
$q : (u, v) \mapsto u|_{\Delta} - v|_{\Delta}.$ From this exact sequence we get 
the following two exact sequences:
\begin{equation}\label{eqn2}
 0 \rightarrow p_1^*(E) \otimes \mathcal{O}_{C \times C}(- \Delta) 
 \rightarrow \pi^*(\mathcal{F}_2(E)) \rightarrow p_2^*(E) 
 \rightarrow 0,
\end{equation}
and
\begin{equation}\label{eqn3}
 0 \rightarrow p_2^*(E) \otimes \mathcal{O}_{C \times C}(- \Delta) 
 \rightarrow \pi^*(\mathcal{F}_2(E)) \rightarrow p_1^*(E) 
 \rightarrow 0
\end{equation}
\cite[Section\, 3]{BN}.

Note that, the action of $\Z/2\Z$ on $C\times C$ lifts to an action on 
$p_1^*(E) \oplus p_2^*(E)$: the fibers of $p_1^*(E) \oplus p_2^*(E)$ are permuted 
in the same way as that of any element of $C \times C.$ Also the bundle 
$\pi^*(\E)$, being pull back from $S^2(C)$, is $\Z/2\Z$-equivariant. And the 
inclusion $\pi^*(\E) \hookrightarrow p_1^*(E) \oplus p_2^*(E)$ is also 
$\Z/2\Z$-equivariant.

\subsection{Proof of Theorem \ref{MT}:}

 $(i)$ By above, it is sufficient to prove the (semi-)stability of $\pi^*(\E)$ on 
 $C \times C$ with respect to the ample divisor $H=[x\times C+C\times x]$. Assume that 
 $E$ is semi-stable, $d\geq r$ and $\pi^*(\E)$ is not semi-stable. Let $A$ be the maximal 
 destabilizing subsheaf of $\pi^*(\E)$. Then, by (\ref{eqn0'}), $\mu_H(A) > \mu_H(\pi^*(\E))\geq 0$. 
 Also $A$ is reflexive and hence locally free. Since $H$ is $\Z/2\Z$-invariant, 
 $A$ has a $\Z/2\Z$-equivariant structure such that the inclusion $A \hookrightarrow 
 \pi^*(\E)$ is $\Z/2\Z$-equivariant (for a proof see \cite[Proposition 4.2.2]{M}).
 Thus the composition of inclusions 
 \begin{equation}\label{eqn4'}
 \psi: A \hookrightarrow \pi^*(\E) \hookrightarrow p_1^*(E) \oplus p_2^*(E)
 \end{equation}
 is also $\Z/2\Z$-equivariant. Hence it is sufficient to show that, for any 
 subsheaf $A$ of $\pi^*(\E)$ with a $\Z/2\Z$ equivariant structure such that the map 
 $\psi: A \hookrightarrow p_1^*(E) \oplus p_2^*(E)$ is $\Z/2\Z$ equivariant, we have 
 $\mu_H(A) \leq \mu_H(\pi^*(\E))$. 
 
 Let $A$ be such a subsheaf of $\pi^*(\E)$ of rank $s$.
 Let $\text{pr}_i: p_1^*(E) \oplus p_2^*(E) \twoheadrightarrow p_i^*(E)$ be the $i$-th 
 coordinate projection and $\psi_i = \text{pr}_i \circ \psi, i= 1, 2.$ 
 Assume that $\text{rank}(\text{Ker}(\psi_1))=s_1$, and  
 $\text{rank}(\text{Im}(\psi_1))=s_2.$ Using $\Z/2\Z$-equivariance of $\psi$, 
 we see that $\text{rank}(\text{Ker}(\psi_2))=s_1$ and $\text{rank}(\text{Im}(\psi_2))=s_2.$ 
 Also notice that $\text{Ker}(\psi_1) \subseteq \text{Im}(\psi_2)$ and $\text{Ker}(\psi_2)
 \subseteq \text{Im}(\psi_1)$ so that $s_1 \leq s_2.$
 
 
 Assume $s_1>0$. Restricting the exact sequence (\ref{eqn2}) to the curves 
 of the form $C \times x$ we get an exact sequence
 $$
 0 \to E \otimes \Oh_{C}(-x) \to \pi^*(\E)|_{C \times x} \to\Oh_C^{\oplus r} \to0.
 $$
 Let $A_2$ (respectively, $A_1$) be the image (respectively, kernel) of the induced 
 map $A|_{C \times x} \to \Oh_C^{\oplus r}.$ Then we have the following commutative 
 diagram:
 \begin{center}
 \begin{tikzpicture}
  \matrix(m)[matrix of math nodes,row sep=3em, column sep=2.5em,text height=1.5ex, 
 text depth=0.25ex]
 {0&E \otimes \Oh_{C}(-x)&\pi^*(\E)|_{C \times x}&\Oh_C^{\oplus r}&0\\
 0&A_1&A|_{C \times x}&A_2&0\\};
 \path[->](m-1-1) edge (m-1-2);
 \path[->](m-1-2) edge (m-1-3);
 \path[->](m-1-3) edge (m-1-4);
 \path[->](m-1-4) edge (m-1-5);
 \path[->](m-2-1) edge (m-2-2);
 \path[->](m-2-2) edge (m-2-3);
 \path[->](m-2-3) edge (m-2-4);
 \path[->](m-2-4) edge (m-2-5);
 \path[->](m-2-2) edge (m-1-2);
 \path[->](m-2-3) edge (m-1-3);
 \path[->](m-2-4) edge (m-1-4);
 \end{tikzpicture}

 \end{center}
 where the rows are exact and the vertical arrows are injections. Note that,
 $\text{rank}(A_1)=s_1$ and $\text{rank}(A_2)=s_2.$ Since $E$ is semi-stable, 
 $\mu(A_1) \leq \frac{d-r}{r}$ and $\mu(A_2) \leq 0.$ Combining these two we get 
 that $\text{deg}(A|_{C \times x}) \leq \frac{s_1(d-r)}{r}$. Similarly, 
 restricting the exact sequence (\ref{eqn3}) to the curves of the form $x\times C$, 
 we get that $\text{deg}(A|_{x \times C}) \leq \frac{s_1(d-r)}{r}$. Thus 
 $\text{deg}_H(A) \leq \frac{2s_1(d-r)}{r}$. Now $s_1 \leq s_2$ implies that 
 $\text{deg}_H(A) \leq \frac{s(d-r)}{r}$, i.e. $\mu_H(A) \leq \mu_H(\pi^*(\E))$. 
 Note that, from the above argument it follows that the inequality for $\text{deg}(A|_{C \times x})$ 
 and $\text{deg}(A|_{x \times C})$ also hold when $s_1=0$. So in this case also, we get 
 $\mu_H(A) \leq \mu_H(\pi^*(\E))$.
 
 \vspace{0.5 cm}
 
 $(ii)$ Let $E$ be slope stable, and $d>r$. Let $\tilde{A}$ be a proper subsheaf 
 of $\E$. Taking double dual, if necessary, we may assume that $\tilde{A}$ is 
 locally free. Then $A:= \pi^*(\tilde{A})$ is a locally free proper subsheaf of 
 $\pi^*(\E)$ such that the inclusion $A \hookrightarrow \pi^*(\E)$ is 
 $\Z/2\Z$-equivariant. Now we proceed as above, and using that $E$ is slope stable and $d>r$ we 
 get $\mu_H(A)<\mu_H(\pi^*(\E))$. Consequently, $\mu_{\{x+C\}}(\tilde{A})<
 \mu_{\{x+C\}}(\E)$.

 This completes the proof.

\begin{Remark}
 Let $\mathcal{M}_1 := M_{_C}(r,d)^s$ denotes the moduli space of stable bundles 
 on $C$ of rank $r$ and degree $d,$ and let 
 $\mathcal{M}_2 := M_{_{S^2(C)}}(c_1, c_2, H, 2r)^s$ denotes the 
 moduli space of stable bundles 
 over $S^2(C)$, where the Chern classes $c_1, c_2$ depend on $d$ and $r$ 
 and can be computed via equation (\ref{eqn0''}), and $H=x+C$. Assume $d>r.$ Then 
 by the above theorem we get that the moduli space $\mathcal{M}_2$ is 
 non-empty and we have a morphism 
 $$
 \Phi : \mathcal{M}_1 \to \mathcal{M}_2, \, \, \, E \mapsto \E.
 $$
 Now using \cite[Theorem 2.1]{BP}, we see that the map $\Phi$ is injective, 
 and the differential map $d\Phi$ is also injective. In other words, $\Phi$ 
 is an immersion.
\end{Remark}

 {\it Acknowledgement:} We would like to thank Prof. D.S. Nagaraj, 
 Prof. A. J. Parameswaran, Prof. T. Ramadas, Prof. V. Balaji, Dr. Sarbeswar Pal 
 for their help and encouragement during this work. We thank the reviewer for useful 
 comments and suggestions to make the exposition better. The second named author was 
 partially supported by NBHM Post-Doctoral Fellowship, DAE (Govt. of India).


\begin{thebibliography}{1111}

\bibitem{ACGH} Arbarello, E.; Cornalba, M.; Griffiths, P. A.; Harris, J.: 
{\it Geometry of algebraic curves.} Vol. I. Grundlehren der Mathematischen 
Wissenschaften, 267. Springer-Verlag, New York, 1985.

\bibitem{BL} Biswas, I.; Laytimi, F.: {\it Direct image and parabolic 
structure on symmetric product of curves.} Jour. Geom. Phys. 61 (2011), 773-780.

\bibitem{BN}  Biswas, Indranil; Nagaraj, D. S.: {\it Reconstructing vector 
bundles on curves from their direct 
image on symmetric powers}. Arch. Math. (Basel) 99 (2012), no. 4, 327-331.

\bibitem{BN1} Biswas, Indranil; Nagaraj, D. S.: {\it Stability of secant bundles 
on second symmetric power of a curve.} Commutative algebra and algebraic geometry 
(CAAG-2010), 13–18, Ramanujan Math. Soc. Lect. Notes Ser., 17, Ramanujan Math. Soc., 
Mysore, 2013.

\bibitem{BP} Biswas, Indranil; Parameswaran, A. J.: {\it Vector bundles on 
symmetric product of a curve.}
J. Ramanujan Math. Soc. 26 (2011), no. 3, 351-355.

\bibitem{DP} Krishanu Dan and Sarbeswar Pal: {\it Semistability of Certain 
Bundles on Second Symmetric Power of a Curve}, Journal of Geometry and Physics, 
103 (2016), 37-42.

\bibitem{HL} Huybrechts, Daniel; Lehn, Manfred: { \it The Geometry of Moduli 
Spaces of Sheaves}. Second Edition. Cambridge University Press. 2010.

\bibitem{LMN} El Mazouni, A.; Laytimi, F.; Nagaraj, D. S.: {\it Secant bundles 
on second symmetric power of a curve}. J. Ramanujan Math. Soc. 26 (2011), no. 2, 
181-194.

\bibitem{M} Ernesto C. Mistretta: {\it Some constructions around stability
of vector bundles on projective varieties } (Ph. D. Thesis) 
\url{http://www.math.unipd.it/~ernesto/}

 \bibitem{Schwarzenberger} Schwarzenberger, R. L. E.: {\it The secant bundle of 
 a projective variety}, Proceedings of the London Mathematical Society, vol. 14 (1964), 
 pp. 369-384.
\end{thebibliography}
\end{document}